\documentclass[11pt,reqno]{amsart}
\usepackage{amsmath}
\usepackage{amsthm}
\usepackage{amssymb}
\usepackage{enumerate}
\usepackage{amscd}
\usepackage{color}
\usepackage{pb-diagram}
\usepackage{graphicx}

\theoremstyle{plain}

\newtheorem{lemma}{Lemma}[section]
\newtheorem{prop}[lemma]{Proposition}

\newtheorem*{prop*}{Proposition}
\theoremstyle{definition}

\theoremstyle{remark}

\numberwithin{equation}{section}

\newcommand{\tr}{\textup{tr}}

\newcommand{\vs}{\vspace{0.2cm}}
\newcommand{\vvs}{\vspace{0.4cm}}

\def \a{\alpha}

\begin{document}
\title[nonnegative quadratic bisectional curvature]{An example of compact K\"ahler manifold with nonnegative quadratic bisectional curvature}
\author{Qun Li}
\address{Department of Mathematics and Statistics\\
                    Wright State University\\
                    3640 Colonel Glenn Highway, Dayton, OH 45435.}
\email{qun.li@wright.edu}
\author{Damin Wu}
\address{Department of Mathematics \\
                The Ohio State University \\
                1179 University Drive, Newark, OH 43055.}
\email{dwu@math.ohio-state.edu}
\author{Fangyang zheng}
\address{Department of Mathematics\\
                    The Ohio State University\\
                    231 West 18th Avenue, Columbus, OH 43210, and\\}
\address{Center for Mathematical Sciences\\Zhejiang University \\
Hangzhou, 310027 China}
\email{zheng@math.ohio-state.edu}
\begin{abstract}
We construct a compact K\"ahler manifold of nonnegative quadratic bisectional curvature, which does not admit any K\"ahler metric of nonnegative orthogonal bisectional curvature. The manifold is a 7-dimensional K\"ahler $C$-space with second Betti number equal to 1, and its canonical metric is a K\"ahler-Einstein metric of positive scalar curvature.
\end{abstract}
\maketitle

\section{Introduction}

\vvs

 In recent years, the condition {\em nonnegative quadratic
bisectional curvature} (which we will denote by $QB\geq 0$) has
drawn more and more attentions; see for example \cite{WYZ}, \cite{FWW},
\cite{CT}, \cite{Zhang}, and \cite{LWZ}. At a point $p$ on a
K\"ahler manifold $(M^n,g)$, this condition is defined by
$$ \sum_{i,j=1}^n R_{i\bar{i}j\bar{j}} \ (x_i-x_j)^2 \geq
0,$$ for any unitary tangent frame $\{ e_1, \ldots , e_n\}$ at $p$
and any real numbers $x_1,$ $\ldots ,$ $x_n$.

\vs

 Note that when the bisectional curvature is nonnegative (denoted
as $B\geq 0$ from now on), namely, when
$R_{X\bar{X}Y\bar{Y}} \geq 0$ for any two type  $(1,0)$
tangent vectors $X,Y$ at $p$, then $QB\geq 0$ at $p$.

\vs

 A condition slightly weaker than $B\geq 0$ is the so-called {\em
nonnegative orthogonal bisectional curvature,} denoted as $B^{\perp }\geq
0$, which requires that $R_{X\bar{X}Y\bar{Y}} \geq 0$ for
any two type $(1,0)$ tangent vectors $X,Y$ at $p$ which satisfy
$X\perp Y$. Clearly, $B^{\perp }\geq 0$ already implies $QB\geq 0$, since
the diagonal terms in the summation vanish, and when $n=2$, $B^{\perp}\geq
0$ and $QB\geq 0$ coincide. However, when $n\geq 3$, $QB\geq 0$ does
not have to have all the orthogonal bisectional curvature terms to be
nonnegative, thus it is weaker than $B^{\perp }\geq 0$, at least from the
algebraic point of view. It is however a totally different question
 whether there will be any compact K\"ahler manifold $(M^n,g)$
($n\geq 3$) with $QB\geq 0$ everywhere such that $M^n$ does not
admit any K\"ahler metric with $B^{\perp }\geq 0$ everywhere.

\vs

 The purpose of this note is exactly to demonstrate the existence
of such a manifold.

\vs

 This topic is of course closely related to the
Frankel-Hartshorne conjectures, which attempts to understand the
elliptic end of the high dimensional uniformization theory. The
famous solution of Mok \cite{Mok} to the generalized Frankel
conjecture states that any compact simply-connected K\"ahler
manifold with $B\geq 0$ everywhere must be biholomorphic to a
compact Hermitian symmetric space. Recently, using the Ricci flow
technique and earlier work of X. X. Chen~\cite{Chen} and Brendle-Schoen~\cite{BS, BS1}, Gu-Zhang~\cite{GZ} proved the following
result:

\vvs

 \noindent {\bf Theorem (Gu-Zhang).} {\em Let $(M^n,g)$ be a
simply-connected
 compact K\"ahler manifold with $B^{\perp }\geq 0$ everywhere. Then $M$ is
 biholomorphic to a compact Hermitian symmetric space. }

\vvs

In other words, the condition $B^{\perp }\geq 0$, although algebraically
weaker than $B\geq 0$, do not generate any new examples, since for
any compact Hermitian symmetric space, its canonical K\"ahler metric
has $B\geq 0$ everywhere.

\vs

 Here we avoided the discussion of non-simply connected cases,
since splitting theorems are already known, in the generalized
Frankel case by the classic slitting theorem of Howard-Smyth-Wu
\cite{HSW}, \cite{W},  and in the generalized Hartshorne case by
Demailly-Peternell-Schneider~\cite{DPS}.

\vs

The generalized Hartshorne conjecture seeks to understand all Fano
manifolds with numerically effective tangent bundles. This class
includes all the K\"ahler C-spaces, namely, all the compact
simply-connected homogeneous K\"ahler manifolds. The conjecture, in
its narrowest sense, states that any compact simply-connected
K\"ahlerian manifold $M^n$ with numerically effective tangent bundle
must be biholomorphic to a K\"ahler C-space (see, for example, \cite[p. 218]{Zheng}).

\vs

Note that for a K\"ahler C-space $M^n$, its canonical
K\"ahler-Einstein metric (which is unique up to a constant multiple)
has $B\geq 0$ everywhere if and only if $M^n$ is Hermitian
symmetric. When $M^n$ is not Hermitian symmetric, any K\"ahler
metric on $M^n$ cannot have $B\geq 0$ everywhere by Mok's Theorem.
In fact, it cannot have $B^{\perp }\geq 0$ everywhere by the recent
theorem of Gu and Zhang. So one has to tolerate some mild negativity
of  bisectional curvature in the quest of generalized Hartshorne
conjecture, at least via differential geometric approach. In light
of this, the condition $QB\geq 0$ comes into play, and it is natural
to ask whether this condition will give a good differential
geometric description of K\"ahler C-spaces.

\vs

By the splitting results of H. Wu et. al., we know that under the
condition $QB\geq 0$ everywhere, any harmonic $(1,1)$ form must be
parallel, thus $M^n$ will admit de Rham decomposition if the second
betti number $b_2>1$. For this reason we should restrict ourselves
to K\"ahler C-spaces with $b_2=1$. This class includes all
irreducible compact Hermitian symmetric spaces.

\vs

Let $M^n$ be an $n$-dimensional K\"ahler $C$-space with $b_2 = 1$. The homogeneous K\"ahler metric on $M^n$ is unique up to a constant multiple, and it is K\"ahler-Einstein. Let us denote this metric by $g_0$.
Let ${\mathcal C}$ be the set of all K\"ahler C-spaces
with $b_2=1$ excluding all irreducible compact Hermitian symmetric
spaces. 
We raise the following

\vvs

\noindent {\bf Conjecture.}

{\em 1). Any $(M^n,g_0)$ in ${\mathcal C}$ has $QB\geq 0$
everywhere.

2). If $(M^n,g)$ is a compact simply-connected  and locally
irreducible K\"ahler manifold with $QB\geq 0$ everywhere, then $M^n$
is biholomorphic to a K\"ahler C-space with $b_2=1$.

3). In 2), $g$ is actually isometric to (a constant multiple of)
$g_0$, if the manifold $M^n$ is not $\mathbb{P}^n$.}

\vvs

In this note, we give an explicit calculation of one special example
$(M^7,g_0)$ in ${\mathcal C}$, and show that it indeed has $QB\geq
0$, thus confirming that the condition is genuinely more tolerant
than $B\geq 0$ or $B^{\perp }\geq 0$. Our computation is brute-force, due
to our lack of knowledge in algebra. We suspect that a more
representation-theoretic computation would establish 1).

\vvs

\noindent {\bf Main Theorem.} {\em The $7$-dimensional K\"ahler
C-space $(B_3, \alpha_2)$ has nonnegative quadratic bisectional
curvature.}

\vvs

\vvs

\section{The curvature of K\"ahler C-spaces with $b_2=1$}

\vvs

K\"ahler C-spaces were studied by H. C. Wang~\cite{Wang}, Borel~\cite{B, B1}, Borel-Hirzebruch~\cite{BH},
and others. A detailed study of the
curvature tensor for such spaces were given by M. Itoh \cite{Itoh}.

Let $G$ be a simply-connected, simple complex Lie
group, and ${\mathfrak g}$ its Lie algebra. Fix a Cartan subalgebra
${\mathfrak h}$ of ${\mathfrak g}$. Let $l=\dim_{\mathbb
C}{\mathfrak h}$, and let $\Delta $ be the root system of ${\mathfrak g}$
with respect to ${\mathfrak h}$. Fix a fundamental root system $\{
\alpha_1, \ldots , \alpha_l\}$ of $\Delta$. It determines an
ordering of the root system $\Delta = \Delta^+ \cup \Delta^-$. We
have
$$ \mathfrak g = \mathfrak h \oplus \bigoplus_{\alpha \in \Delta } \mathfrak g_{\alpha },$$
where $\mathfrak g_{\alpha }$ is the root space corresponding to
$\alpha $, satisfying $ [\mathfrak g_{\alpha }, \mathfrak
g_{\beta } ] \subseteq \mathfrak g_{\alpha +\beta }$ for any two
roots $\alpha$, $\beta \in \Delta$.

\vs

Now let us fix an integer $r$ with $1\leq r\leq l$. Denote by
$$\Delta^+_r(k)=\{ \sum_i n_i \alpha_i  \in \Delta^+ | \ n_r=k\}, \ \
\ \ \ \ \Delta^+_r = \bigcup_{k>0} \Delta^+_r(k) .$$ Let $P\subseteq
G$ be the subgroup whose Lie algebra is
$$ \mathfrak p = \mathfrak h \oplus \bigoplus_{\alpha \in \Delta
\setminus \Delta^+_r} \mathfrak g_{\alpha } , $$ then $P$ is a
parabolic subgroup, namely, the complex manifold $M=G/P$ is compact.
This gives a K\"ahler C-space with $b_2=1$. Conversely, any K\"ahler
C-space with $b_2=1$ is given this way. We will denote this space by
$(\mathfrak g, \alpha_r)$.

\vs

Table 1 on page 55 of \cite{Itoh} gives the list of all K\"ahler
C-spaces with $b_2=1$, using the Dynkin diagrams. The double circled
ones are the Hermitian symmetric ones. From the table, we see that the
simplest non-symmetric example would be $M^7=(B_3, \alpha_2)$, which
will be our example in this note.

\vs

In the remainder of this section, we will follow Itoh's notations
and calculations \cite{Itoh}, and collect the necessary formulas that we need later for computing the
curvature of K\"ahler C-spaces with $b_2=1$. Let $M^n=(\mathfrak g, \alpha_r)$ be a K\"ahler C-space
with $b_2=1$. Denote 
$$ {\mathfrak m}^{+}_k = \bigoplus_{\alpha \in \Delta^+_r(k)} \mathfrak
g_{\alpha }, \ \ \ \ \ \ \mathfrak m^+ = \bigoplus_{k\geq 1}
{\mathfrak m}^{+}_k, \ \ \ \ \ \ $$ and define $\mathfrak m^-_k$ and
${\mathfrak m}^-$ similarly. Also, let $\mathfrak t= \mathfrak h
\oplus \bigoplus_{\alpha \in \Delta^+_r(0)} (\mathfrak g_{\alpha }
\oplus {\mathfrak g}_{-\alpha } )$. Then ${\mathfrak g}=\mathfrak t
\oplus \mathfrak m$, where ${\mathfrak m}=\mathfrak m^+ \oplus
\mathfrak m^-$, and
$$ [\mathfrak t,\mathfrak m^{\pm}_k]\subseteq \mathfrak m^{\pm }_k , \ \
\ \ [\mathfrak m^{\pm}_k,\mathfrak m^{\pm}_l]\subseteq \mathfrak
m^{\pm}_{k+l}, \ \ \ \ [\mathfrak m^+_k,\mathfrak m^-_k]\subseteq
\mathfrak t, $$ and for any $k>l>0$, it holds \ $[\mathfrak
m^+_k,\mathfrak m^-_l]\subseteq \mathfrak m^+_{k-l}$, \ $[\mathfrak
m^-_k,\mathfrak m^+_l]\subseteq \mathfrak m^-_{k-l}$. The space
$\mathfrak m^+$ can be identified with the holomorphic tangent space
of $M^n$. Let $K$ be the Killing form of ${\mathfrak g}$. Let $E_{\alpha }$ be a
Weyl canonical basis of ${\mathfrak g}$, namely,
$E_{\alpha } \in {\mathfrak g}_{\alpha }$ for each $\alpha \in \Delta$, and
\begin{equation} \label{eq:test1}
K(E_{\alpha }, E_{-\alpha })=-1, \qquad  N_{\alpha , \beta } = N_{-\alpha , -\beta },
\end{equation}
for any $\alpha, \beta \in \Delta^+$, where $N_{\alpha , \beta }$ is determined by
$[E_{\alpha }, E_{\beta }]=N_{\alpha , \beta } E_{\alpha +\beta }$.

\vs

The canonical metric $g_0=\langle \ , \ \rangle $ is given by
\begin{equation} \label{eq:test2}
 g_0(X,\bar{Y})=\langle X, \bar{Y} \rangle = - k K(X,\bar{Y})
 \end{equation}
for any $X$, $Y\in \mathfrak m^+_k$, where the complex conjugation on $M$ is determined by
$\bar{E_{\alpha }}=E_{-\alpha }$ for each $\alpha$.
In the following, we will assume that
$$ X\in \mathfrak m^{+}_i, \ \ Y \in  \mathfrak m^{+}_j, \ \ Z\in  \mathfrak
m^{+}_k, \ \ W \in \mathfrak m^{+}_l,$$ where $i,j,k,l$ are any positive integers,
and compute the curvature component $R_{X\bar{Y}Z\bar{W}}$ of the canonical
metric $g_0$. From \cite{Itoh}, we have
$$ R(X,\bar{Y})Z= [\Lambda (X), \Lambda (\bar{Y})]Z-\Lambda ([X,\bar{Y}]_{\mathfrak m})Z-[[X,\bar{Y}]_{\mathfrak t},Z], $$
where
$$ \Lambda (X)Y=\frac{j}{i+j}[X,Y], \ \ \ \ \ \Lambda(\bar{X})Y=[\bar{X},Y]_{\mathfrak m^{ +}}. $$
Like in \cite{Itoh}, a straight forward computation yields the following formulas.

\begin{prop} \label{pr:2.1}
On a K\"ahler C-space $(M^n,g_0)$ with $b_2=1$, then for any
$$ X\in \mathfrak m^{+}_i, \ \ Y \in  \mathfrak m^{+}_j, \ \ Z\in  \mathfrak
m^{+}_k, \ \ W \in \mathfrak m^{+}_l,$$
the curvature components are given by
\begin{equation} \label{eq:RXYZW}
 \begin{split}
		R_{X\bar{Y}Z\bar{W}} = \ & (k-j)\xi_{k-j}K([X,\bar{W}], [\bar{Y}, Z]) - \frac{kl}{i+k}K([X,Z], [\bar{Y}, \bar{W}]) \\
& +  (k\xi_{i-j}+l\xi_{j-i}+l\delta_{ij}\delta_{kl})K([X,\bar{Y}],[Z,\bar{W}]),
 \end{split}
\end{equation}
if $i+k=j+l$, and $R_{X\bar{Y}Z\bar{W}} = 0$ otherwise. Here $\xi_q=1$ for $q>0$ and $\xi_q=0$ for $q\leq 0$. \qed
\end{prop}

Note that by the invariance of $K$ and the Jacobi identity, we see that the first (or the third) term on the right hand side of \eqref{eq:RXYZW}
can be expressed as a linear combination of the other two terms.

\vs

In the following, we will assume that $\mathfrak g$ satisfies the condition
\begin{equation} \label{eq:contact}
 \Delta^+_r(k)=\emptyset , \ \ \textup{for all $k\geq 3$}.
 \end{equation}
This condition is satisfied by all four classical sequences $A$, $B$, $C$, $D$ for all $r$, and for some of the exceptional cases.
Note that the case ${\mathfrak m}^+=\mathfrak m^+_1$ corresponds exactly to all the irreducible Hermitian symmetric cases.

\vs

For the sake of simplicity, we will assume that $M^n=(\mathfrak g, \alpha_r)$ is of {\em contact type,} namely, the condition \eqref{eq:contact} holds
and $\Delta^+_r(2)$ consists of only one element.  Applying Proposition~\ref{pr:2.1} to
this contact case, we get

\begin{prop} \label{pr:2.2}
Let $M^n$ be a K\"ahler C-space with $b_2=1$ of the contact type, then for any $X,Y,Z,W$
in ${\mathfrak m}^+_1$, and $U\in \mathfrak m^+_2$, we have
\begin{align*}
R_{U\bar{U}U\bar{U}} & = 2K([U, \bar{U}], [U, \bar{U}]),\\
R_{X\bar{Y}U\bar{U}} & = - K([X, \bar{U}],[\bar{Y},U]) = K([X, \bar{Y}],[U, \bar{U}]),\\
R_{X\bar{Y}Z\bar{U}} & = R_{X\bar{U}U\bar{U}} \ = R_{X\bar{U}Z\bar{U}} \ = 0,\\
R_{X\bar{Y}Z\bar{W}} & = -\frac{1}{2} K([X, Z],[\bar{Y},\bar{W}]) + K([X, \bar{Y}],[Z, \bar{W}])\\
& = \frac{1}{2} K([X, Z],[\bar{Y},\bar{W}]) +  K([X, \bar{W}],[Z, \bar{Y}]).
\end{align*} \qed
\end{prop}

\section{The curvature tensor of $(B_3,\alpha_2)$}

\vvs

Let us now consider the specific case of $M^7=(B_3, \alpha_2)$, where $B_3 = \mathfrak{so}_7(\mathbb{C})$ is the Lie algebra of
the Lie group $SO_{7}(\mathbb{C})$. We will regard $B_3$ as the space of all $7\times 7$ complex matrices $A=(a_{ij})$ such that
$a_{ij}=-a_{j'i'}$ for all
$1\leq i,j\leq 7$, where here and from now on we will write $i'=8-i$. A Cartan subalgebra is given by all the diagonal matrices
\[
   \mathfrak{h} = \big\{\textup{diag}(a_1, a_2, a_3, 0, - a_3, - a_2, - a_1)\mid a_i \in \mathbb{C}, 1 \le i \le 3\big\}.
\]
The root system with respect to this $\mathfrak h$ is
\[
   \Delta = \{ \pm \varepsilon_i \pm \varepsilon_j \mid 1 \le i,j \le 3, i \ne j \} \cup \{ \pm \varepsilon_i \mid 1 \le i \le 3 \},
\]
where $\varepsilon_i \in \mathfrak{h}^*$ is defined by
\[
   \varepsilon_i\big(\textup{diag}\{b_1, b_2, b_3, 0, -b_3, -b_2, -b_1\}\big) = b_i, \qquad 1 \le i \le 3.
\]
Note that $\mathfrak{g} = \mathfrak{so}_7(\mathbb{C})$ has three simple roots $\{\alpha_1, \alpha_2, \alpha_3 \}$, which are given by
\[
   \a_1 = \varepsilon_1 - \varepsilon_2, \quad \a_2 = \varepsilon_2 - \varepsilon_3, \quad \a_3 = \varepsilon_3.
\]
Let $\Delta^+$ and $\Delta^-$ be the subsets of $\Delta$ consisting of positive and negative roots, respectively. We have
\begin{align*}
\Delta_2^+(0) & = \{ \alpha_1, \alpha_3 \} \ = \ \{ \varepsilon_1-\varepsilon_2, \ \varepsilon_3  \} ,  \\
   \Delta_2^+ (1) & = \{ \alpha_2, \alpha_2 + \alpha_1, \alpha_2 + \alpha_3, \alpha_2 + 2\alpha_3, \alpha_2 + \alpha_1 + \alpha_3,
   \alpha_2 + \alpha_1 +
   2\alpha_3 \}  \\
   & = \{\varepsilon_2 - \varepsilon_3, \varepsilon_1 - \varepsilon_3, \varepsilon_2, \varepsilon_2 + \varepsilon_3, \varepsilon_1,
   \varepsilon_1 + \varepsilon_3\},  \\
   \Delta_2^+(2) & =\{2\alpha_2 + \alpha_1 + 2\alpha_3 \} = \{ \varepsilon_1 + \varepsilon_2 \},  \\
   \Delta_2^+(k) & = \emptyset, \quad \textup{for\ $k \ge 3$.}
\end{align*}
Therefore, $\varDelta_2^+ = \varDelta_2^+(1) \sqcup \varDelta_2^+(2)$ and $M^7$ is of contact type.

\vs

For convenience, let us denote by $e_{ij}$ the square matrix whose $(i,j)$-entry is $1$ and other entries are zero. Let
\[
   F_{ij} = e_{ij} - e_{j'i'}, \qquad 1 \le i,j \le 7.
\]
Then, $\mathfrak{h}$ is spanned by $\{F_{ii}\mid 1 \le i \le 3\}$ over $\mathbb{C}$. Furthermore, the root vectors  can
be explicitly given by
\begin{align*}
    &\mathfrak{g}_{\varepsilon_i-\varepsilon_j} = \mathbb{C}F_{ij}, \quad \mathfrak{g}_{\varepsilon_i+\varepsilon_j} = \mathbb{C}
    F_{ij'}, \quad \mathfrak{g}_{-\varepsilon_i-\varepsilon_j} = \mathbb{C} F_{i'j}, \\
    &\mathfrak{g}_{\varepsilon_i} = \mathbb{C} F_{i\,4}, \quad \mathfrak{g}_{-\varepsilon_i} = \mathbb{C}F_{4\,i}, \qquad
    \textup{for all $1 \le i, j \le 3$}.
 \end{align*}
We compute the trace form
\[
  (F_{ij}, F_{kl})_{\mathbb{C}^7} = \tr_{\mathbb{C}^7}(F_{ij}F_{kl}) = 2(\delta_{il}\delta_{jk} - \delta_{ik'}\delta_{jl'}).
\]
Since $\mathfrak{so}_7(\mathbb{C})$ is a simple Lie algebra, any two invariant symmetric bilinear form on $\mathfrak{so}_7(\mathbb{C})$
are proportional. Thus, by rescaling some constant, we can assume that the Killing form on $\mathfrak{so}_7(\mathbb{C})$ satisfies
\begin{equation} \label{eq:Killing}
   K(F_{ij},F_{kl}) =  \delta_{il}\delta_{jk} -  \delta_{jl'}\delta_{ik'}.
\end{equation}
Moreover, observe that
\begin{equation} \label{eq:bkt}
   [F_{ij}, F_{kl}] = \delta_{jk}F_{il} -\delta_{il}F_{kj} -\delta_{j'l}F_{ik'} - \delta_{ik'}F_{j'l}.
\end{equation}
We will form a unitary tangent frame $\{E_1, \ldots , E_6; E_7\}$ (in the order) by  
\begin{equation*} 
 \{ F_{13}, F_{14}, F_{15}, F_{23}, F_{24}, F_{25}; \frac{1}{\sqrt{2}}F_{16} \} 
\end{equation*}
, and let their complex conjugation $\{ \bar{E}_1, \ldots , \bar{E}_6; \bar{E}_7\}$ be
$$ \{ -F_{31}, -F_{41}, -F_{51}, -F_{32}, -F_{42}, -F_{52}; -\frac{1}{\sqrt{2}}F_{61} \} .$$

In the following, for the benefit of calculations in the next section, we will let $a, b, c, d$ be the indices from $\{ 1, 2\}$ and $i,j,k,l$ the indices from $\{ 3,4,5\}$. We will also let $X=F_{ai}$, $Y=F_{bj}$, $Z=F_{ck}$, $W=F_{dl}$ be vectors in $\mathfrak m^+_1$ and
$U=\frac{1}{\sqrt{2}}
F_{16}$ be
in $\mathfrak m^+_2$. Note that $F_{ij}=-F_{j'i'}$ for any $i,j$. By \eqref{eq:bkt}, we get
\begin{align*}
[X, \bar{Y}] & = \delta_{ab}F_{ji} - \delta_{ij} F_{ab}, \qquad
[X,Z]  = -\delta_{ik'} F_{ac'},  \\
[Z,\bar{W}] & = \delta_{cd}F_{lk} - \delta_{kl}F_{cd}, \qquad
[\bar{Y}, \bar{W}]  = - \delta_{jl'} F_{b'd}, \\
[U,\bar{U}] & = - \frac{1}{2}(F_{22}+ F_{11}) .
\end{align*}
So by Proposition~\ref{pr:2.2} and formula \eqref{eq:Killing}, we get from a straight forward computation that
\begin{align}
R_{U\bar{U}U\bar{U}} & = 2K([U,\bar{U}], [U,\bar{U}]) = 1, \notag \\
R_{X\bar{Y}U\bar{U}} & = \frac{1}{2} \delta_{ij} \delta_{ab}, \notag \\
R_{X\bar{Y}Z\bar{W}} & = -\frac{1}{2}(\delta_{ad}\delta_{bc}+ \delta_{ab}\delta_{cd}) \delta_{ik'}\delta_{jl'}
+ \delta_{ad}\delta_{bc}\delta_{ij}\delta_{kl} + \delta_{ab}\delta_{cd}\delta_{il}\delta_{jk}. \label{eq:3.3}
\end{align}
 The result below follows immediately.
\begin{prop}
The canonical metric $g_0$ on $M^7 = (B_3,\alpha_2)$ is K\"ahler-Einstein, and its Ricci curvature is identically equal to $4$. \qed
\end{prop}

\vvs

\vvs

\section{Nonnegativity of the quadratic bisectional curvature}

\vvs

In this section, we shall show that the $7$-dimensional K\"ahler C-space $(M^7,g_0)$ given by
$(B_3,\alpha_2)$ indeed has nonnegative quadratic bisectional curvature. 
In other words, for any unitary tangent frame
$\{ e_{a}: 1\leq a\leq 7\}$
and any real numbers $x_1,\ldots , x_7$, the quantity $QB : =  \sum_{a,b=1}^7 R(e_a, \bar{e}_a, e_b, \bar{e}_b) (x_a-x_b)^2$ is
always nonnegative.

\vs

Write $e_a=\sum_{i=1}^7 T_{ai}E_i$ for some $T\in U(7)$, where $\{ E_1, \ldots , E_7\}$ is the Weyl basis mentioned in the \S 3.
Since the metric $g_0$ is
Einstein with Ricci curvature equal to $4$,  we have
\begin{align*}
\frac{1}{2}QB & = \sum_a R_{a\bar{a}}x_a^2 - \sum_{a,b} R_{a\bar{a}b\bar{b}}x_ax_b \\
& = 4\sum_a x_a^2 - \sum R_{i\bar{j}k\bar{l}} T_{ai}\bar{T}_{aj}T_{bk}\bar{T}_{bl}x_ax_b \\
& = 4 \sum_{i,j} |P_{i\bar{j}}|^2 - \sum_{i,j,k,l}R_{i\bar{j}k\bar{l}} P_{i\bar{j}} P_{k\bar{l}} \\
& = 4 |P|^2 - \Box.
\end{align*}
where $P=\ ^t\!T \Lambda_x \bar{T}$ is a Hermitian symmetric matrix. Here we wrote $\Lambda_x=\mbox{diag} \{ x_1, \ldots , x_7\}$
and $|P|^2=\sum_{i,j} |P_{i\bar{j}}|^2 = \sum_a x_a^2$. Write
\[
P=\left[ \begin{array}{ll} P' & \xi \\ ^t\!\bar{\xi} & t \end{array} \right].
\]
We have $|P|^2 = |P'|^2+ 2|\xi|^2 + t^2$. Since under the Weyl frame $\{E_i\}$, we have $R_{7\bar{7}7\bar{7}}=1$
and $R_{i\bar{j}7\bar{7}}=\frac{1}{2}\delta_{ij}$ for any $1\leq i,j\leq 6$, we know that
\begin{align*}
\Box & = \sum_{i,j,k,l=1}^6 R_{i\bar{j}k\bar{l}} P_{i\bar{j}} P_{k\bar{l}} + t^2 + t \cdot \mbox{tr}(P') + |\xi |^2\\
& = \Box ' + t^2 + t \cdot \mbox{tr}(P') + |\xi |^2.
\end{align*}
Plugging these into the formula for $QB$, we get
$$ \frac{1}{2} QB = 7|\xi |^2 + 3t^2 - t \cdot \mbox{tr}(P')+ (4|P'|^2-\Box ').$$
So $QB\geq 0$ for any Hermitian matrix $P$ is equivalent to
\begin{equation} \label{eq:4.1}
\Phi := 4|P'|^2 - \Box ' -\frac{1}{12} (\mbox{tr}(P'))^2 \geq 0
\end{equation}
for any $6\times 6$ Hermitian matrix $P'$.

\vs

Writing $\{ E_1, \ldots E_6\}$ as $\{ F_{ai} : 1\leq a\leq 2, \ 3\leq i\leq 5\}$ 
and using formula~\eqref{eq:3.3}, we get
\begin{align*}
\Box ' & = \sum_{a,b=1}^2 \sum_{i,k=3}^5 \big( -\frac{1}{2} P_{ai\overline{ak}} P_{bi'\overline{bk'}}
-\frac{1}{2} P_{ai\overline{bk}} P_{bi'\overline{ak'}} \\
& +  P_{ai\overline{bi}} P_{bk\overline{ak}} + P_{ai\overline{ak}} P_{bk\overline{bi}} \ \big).
\end{align*}
Here as before, $i'=8-i$,  and we use double indices for $P$, e.g., $P_{13\overline{13}} = P_{1\bar{1}}$. Now let $A$, $B$, $C$ be the $3\!\times\!3$ matrices given by
$$ A_{ij}=P_{1i\overline{1j}}, \ \ \ B_{ij}=P_{1i\overline{2j}}, \ \ \ C_{ij}=P_{2i\overline{2j}};$$
then we can write $\Box '$ as
\begin{align*}
&  -\frac{1}{2} \langle A+C, A+C\rangle -\frac{1}{2}\big( \langle A,A\rangle + \langle C,C\rangle + \langle B,B^*\rangle +
\langle B^*,B\rangle \big) \\
& + (\mbox{tr} A)^2+ (\mbox{tr} C)^2 +2|\mbox{tr} B|^2 + |A+C|^2,
\end{align*}
where $|A|^2=\sum_{i,j=3}^5 |A_{ij}|^2$, $\langle A,B\rangle = \sum_{i,j=3}^5 A_{ij}B_{i'j'}= \langle B,A\rangle $, and $B^*= \ ^t\!\bar{B}$.
Since $\mbox{tr}(P')=\mbox{tr}A+\mbox{tr}C$, and
$|P'|^2=|A|^2+|C|^2+2|B|^2$, we get
\begin{align*}
\Phi & = 4(|A|^2+|C|^2+2|B|^2)-\Box ' -\frac{1}{12}( \mbox{tr}A+\mbox{tr}C)^2 \\
& = \Phi_1+ \Phi_2,
\end{align*}
where
\begin{align*}
 \Phi_1   = \ &  4|A|^2 + 4|C|^2 + \langle A, A\rangle + \langle C, C\rangle + \langle A, C\rangle
  - (\mbox{tr}A)^2 \\ & - (\mbox{tr}C)^2 -   |A+C|^2 - \frac{1}{12}(\mbox{tr}A+ \mbox{tr}C)^2,\\
 \Phi_2  =  \ &  8|B|^2 +  \langle B, B^* \rangle -2\ |\mbox{tr}B|^2.
\end{align*}
From
$$ \langle B,B^*\rangle = |B_{44}|^2+|B_{35}|^2+|B_{53}|^2 +2\mbox{Re}(B_{33}\bar{B}_{55}+ B_{34}\bar{B}_{45} + B_{43}\bar{B}_{54}),$$
we see that $\langle B, B^* \rangle \geq - |B|^2$. This together with the fact that $|\mbox{tr}B|^2\leq 3|B|^2$ imply that
$$\Phi_2\geq (8-1-6)\ |B|^2\geq 0.$$
To deal with $\Phi_1$, let us write $\Phi_1=\Phi_1'+\Phi_1''$, where $\Phi_1'$ is the part of $\Phi_1$ coming from the $(34)$, $(43)$, $(45)$,
and $(54)$ entries of $A$ and $C$, and $\Phi_1''$ the rest. Then, we have
\begin{align*}
 \Phi_1'  = & \ 8|A_{34}|^2+8|A_{45}|^2 + 8|C_{34}|^2+ 8|C_{45}|^2 +4\mbox{Re}(A_{34}\bar{A}_{45}+C_{34}\bar{C}_{45}) \\
 &  + 2\mbox{Re}(A_{34}\bar{C}_{45}+ C_{34}\bar{A}_{45})) -2|A_{34}+C_{34}|^2 -2|A_{45}+C_{45}|^2.
\end{align*}
Using Cauchy-Schwartz inequality, it is easy to see that
$$ \Phi_1'\geq (|A_{34}|^2+|A_{45}|^2 + |C_{34}|^2+ |C_{45}|^2) \geq 0. $$
So we are only left with the verification of the nonnegativity of $\Phi_1''$. Let us denote by $a=(a_1,a_2,a_3)$ the three diagonal
elements of $A$, let $c=(c_1,c_2,c_3)$ be the three diagonal elements of $C$, and write $x=A_{35}$, $y=C_{35}$. Then we have
$\Phi_1''=\Phi_{11}+\Phi_{12}$, where
\begin{align*}
\Phi_{11} & = 10|x|^2+10|y|^2+2\mbox{Re}(x\bar{y})-2|x+y|^2, \\
\Phi_{12} & = 4|a|^2+4|c|^2 +(2a_1a_3+a_2^2) + (2c_1c_3+c_2^2) + (a_1c_3+a_3c_1+a_2c_2)\\
&  - t_a^2-t_c^2 -|a+c|^2 -\frac{1}{12}(t_a+t_c)^2,
\end{align*}
where $t_a=a_1+a_2+a_3$ and $|a|^2=a_1^2+a_2^2+a_3^2$, and likewise for $c$. Clearly, $\Phi_{11}\geq 0$, and $\Phi_{12}$ is a
homogeneous polynomial
of degree $2$ in $a_i$ and $c_i$, $i=1,2,3$. It is not hard to see that $\Phi_{12}=(a,c)D\ ^t\!(a,c)$, where $D$ is the real
$6\!\times\!6$ symmetric matrix given by
\[ D=\left[ \begin{array}{ll} G & H \\ H & G\end{array} \right],
\]
in which $G=3I-\frac{13}{12}L+J$, $H=-I-\frac{1}{12}L+\frac{1}{2}J$, with
\[ L=\left[ \begin{array}{lll} 1 & 1 & 1 \\ 1 & 1 & 1 \\ 1& 1 & 1 \end{array} \right], \ \ \ \
J = \left[ \begin{array}{lll} 0 & 0 & 1 \\ 0 & 1 & 0 \\ 1& 0 & 0 \end{array} \right],
\]
and $I=I_3$ the identity matrix. We want to show that $D\geq 0$, or equivalently, the $3\!\times\!3$ matrix
\[ G - H G^{-1}H \geq 0. \]
Consider the matrix $T\in O(3)$ given by
\[ T= \frac{1}{\sqrt{6}} \left[ \begin{array}{ccc} \sqrt{3} & 1 & \sqrt{2} \\ 0 & -2 & \sqrt{2} \\
-\sqrt{3} & 1 & \sqrt{2} \end{array} \right]. \]
We have
\[ ^t\!TLT= 3 \left[ \begin{array}{lll} 0 & 0 & 0 \\ 0 & 0 & 0 \\ 0 & 0 & 1 \end{array} \right] , \ \ \ \
 ^t\!TJT= \left[ \begin{array}{ccc} -1 & 0 & 0 \\ 0 & 1 & 0 \\ 0 & 0 & 1 \end{array} \right].
\]
Therefore, we get
\[ ^t\!TGT =\left[ \begin{array}{lll} 2 & 0 & 0 \\ 0 & 4 & 0 \\ 0 & 0 & \frac{3}{4} \end{array} \right] ,
\ \ \ \ ^t\!THT =\left[ \begin{array}{ccc} -\frac{3}{2} & 0 & 0 \\ 0 & -\frac{1}{2} & 0 \\ 0 & 0 & -\frac{3}{4} \end{array} \right].
\]
Thus,
\[ ^t\!T(G-HG^{-1}H)T =\left[ \begin{array}{lll} 2 & 0 & 0 \\ 0 & 4 & 0 \\ 0 & 0 & \frac{3}{4} \end{array} \right] -
\left[ \begin{array}{ccc} \frac{9}{8} & 0 & 0 \\ 0 & \frac{1}{16} & 0 \\ 0 & 0 & \frac{3}{4} \end{array} \right]
\geq 0.
\]
This establish the nonnegativity of $\Phi_{12}$, hence the entire $\Phi$, and by \eqref{eq:4.1} we have shown that $(B_3,\alpha_2)$ has
nonnegative quadratic bisectional curvature.

\vvs

\vvs

\vvs

\vvs

\vvs

\end{document}